\documentclass[numreferences]{kluwer}
\usepackage{graphicx}
\usepackage{setspace}
\usepackage{enumitem}
\usepackage{amsthm}
\newtheorem{Thm}{Theorem}[section]

\newtheorem{Alg}[Thm]{Algorithm}

\begin{document}
	\begin{opening}
		\title{Trigonometric-Interpolation Based Approach for Second-Order Volterra Integro-Differential Equations}
		\author{Xiaorong \surname{Zou}\email{xiaorong.zou@bofa.com}}
		\institute{Global Market Risk Analytic, Bank of America}
		\runningauthor{X. Zou}
		\runningtitle{TIBA for Second Order VIDE}
		\date{May 6, 2025}
		\classification{MSC2000}{Primary 65T40; Secondary 45B05}
		\keywords {Trigonometric Interpolation, Volterra Integro-differentiable Equation, Fast Fourier Transform}

\begin{abstract}
The	trigonometric interpolation has been recently applied to solve a second order Fredholm integro-differentiable equation (FIDE). It achieves high accuracy with a moderate size of grid points and effectively address singularities of kernel functions. In addition, it work well with general boundary conditions and the framework can be generalized
to work for FIDEs with high order ODE component.  In this paper, we apply the same idea to develop an algorithm for the solution of a second order Volterra integro-differentiable equation (VIDE) with same advantages as in the study of FIDE. The numerical experiments with various boundary conditions are conducted with decent performances as expected.			
\end{abstract}
\end{opening}
\section{Introduction}\label{sec:vide_intro}
A $2$-dim trigonometric interpolation algorithm was recently introduced to approximate non-periodic functions in $2$-dim space and used to develop the solution of Fredholm integro-differentiable equation (FIDE) in \cite{zou_tri_FIDE}. 

In this paper, we follow similar approach to develop an algorithm for the solution of the following second order Volterra integro-differentiable equation (VIDE).  
\begin{equation}\label{eq:vide_ode}
	y''(x)=p(x)y'+q(x)y(x) +r(x) + \mu(x)\int^x_sk(x,t)y(t)dt, \quad x\in [s,e].  
\end{equation}
\begin{eqnarray}
	d_{11}y(s) + d_{12}y'(s) +d_{13}y(e) + d_{14}y'(e) &=& \alpha , \label{eq:vide_bdy1}\\ 
	d_{21}y(s) + d_{22}y'(s) +d_{23}y(e) + d_{24}y'(e) &=& \beta, \label{eq:vide_bdy2} 
\end{eqnarray}
where $p,q,r,\mu$ is continuous differentiable on $[s, e]$, the kernel function $k(x,y)$ is integrable with respect to $y$ on $[s,t]$ for any $x\in [s,e]$, the rank of matrix $D:=(d_{ij})_{1\le i\le 2, 1\le j\le 4}$ is $2$, and $\alpha, \beta$ are two real numbers.  

IDEs appear in mathematical physics, applied mathematics and engineering. In general, analytical solutions of IDEs are not feasible and certain approximations are required. Some numerical solutions in recent years can be found in (\cite{ide_23}-\cite{ide_36}). As in the study of FIDE, the proposed new algorithm converts a VIDE into a linear algebraic system with some advantages.  High accuracy can be achieved with moderate size of grid points as shown by numerical examples conducted in Section \ref{sec:numerical_examples}; it can effectively address singularity issue in kernel function as described in Subsection \ref{subsec:fide_en}; it works well for various boundary conditions.  Finally, it can be enhanced to work for IDEs with high order ODE component. 

The rest of paper is organized as follows.  
Section \ref{sec:fide} is denoted to develop Algorithm \ref{Alg:VIDE} to convert VIDE (\ref{eq:vide_ode}-\ref{eq:vide_bdy2}) into a linear algebraic system. It begins  with reviewing related results in \cite{zou_tri_FIDE} to attack the ODE component in VIDE.  Integral component is addressed in Subsection \ref{subsec:fide} and \ref{subsec:fide_en}  for continuous and integrable kernels respectively.  The proposed algorithm is summarized in Subsection \ref{subsec:transform_linear_algebraic_system}. Section \ref{sec:numerical_examples} demonstrates quite a few numerical examples to assess the performance of Algorithm \ref{Alg:VIDE}. Tests covers four types of boundary conditions with both continuous kernels and integrable kernels with singularity. The results show decent accuracy for all covered scenarios with moderate grid point size. The summary is made on Section \ref{sec:summary}.
\section{Development of the Algorithm}\label{sec:fide}
In this section,  we aim to develop Algoirthm \ref{Alg:VIDE} to solve VIDE (\ref{eq:vide_ode}-\ref{eq:vide_bdy2}).
Let $p(x)$,$q(x)$, $r(x)$, $K(x,t)$ be defined in Eq (\ref{eq:vide_ode}), define
\[
f(x,v,u) = p(x)u+q(x)v +r(x), \quad  g(x,v) = \mu(x)\int^x_sK(x,t)v(t)dt
\]
VIDE (\ref{eq:vide_ode}) is reduced to linear ODE if the integration component $g(x,v)$ vanishes and has been solved in \cite{zou_tri_IV} by transforming ODE and boundary condition (\ref{eq:vide_bdy1}-\ref{eq:vide_bdy2}) into linear algebraic system.  We shall follow the same method and convert VIDE (\ref{eq:vide_ode})  to a linear algebraic system by leveraging what has been developed in \cite{zou_tri_IV} to handle  ODE component $f$ and attack integration component $g$ in Subsection \ref{subsec:fide} and \ref{subsec:fide_en} based on smoothness of kernel function $k$.

To apply trigonometric interpolation for non-periodic functions, 
we assume that $p(x)$, $q(x)$ and $r(x)$ are continuous differentiable on $[s-\delta, e+\delta]$ for certain $\delta>0$ and $K(x,t)$ are continuous differentiable on $[s-\delta, e+\delta]^2$. 
By parallel shifting if needed,  we assume $s=\delta$ without loss of generality. Let $h$  be a cut-off function specified in \cite{zou_tri_FIDE} and 
construct $f_h(x,v,u)$ and $g_h(x,t)$ as follows
\[
f_h(x,v,u) = f(x,v,u)h(x), \quad g_h(x,v) =h(x)g(x,v).
\]
Consider a solution $v(x)$ of the following VIDE system
\begin{eqnarray}
	v''(x) &=& f_h(x, v,v') + g_h(x,v),  \quad x\in [0,b] \label{eq:nonlinear_ode_order2_F},\\
	\alpha  &=& d_{11}v(s) + d_{12}v'(s) +d_{13}v(e) + d_{14}v'(e), \nonumber \label{eq:nonlinear_ode_order2:diri_F} \nonumber\\ 
	\beta  &=& d_{21}v(s) + d_{22}v'(s) +d_{23}v(e) + d_{24}v'(e). \nonumber 
\end{eqnarray}
It is clear that $v(x)|_{[s,e]}$ solves VIDE (\ref{eq:vide_ode}-\ref{eq:vide_bdy2}). Define $u(x):=v'(x)$ and $z(x):=v''(x)$.  By Eq (\ref{eq:nonlinear_ode_order2_F}),  $z(x)$ and its derivatives $z^{(k)}$ vanish at boundary points $\{0,b \}$, hence it can be smoothly extended as an odd periodic function with period $2b$ and be approximated by trigonometric polynomial.  Assume 
\begin{equation}\label{eq:vide_z}
	\tilde{z}_M(x) = \sum_{0\le j<M}b_j \sin \frac{j\pi x}{b} 
\end{equation}
is an interpolant of $z(x)$ with $N$ equispaced grid points over $[-b,b]$. $u$ and $v$ can be derived accordingly 
\begin{eqnarray}
	\tilde{u}_M(x ) &=& a_0 -   \frac{b}{\pi }\sum_{1\le j <M} \frac{b_j}{j} \cos  \frac{j\pi x}{b},  \label{eq:vide_u} \\
	\tilde{v}_M(x ) &=& a_{-1} +  a_0 x -  (\frac{b}{\pi })^2 \sum_{1\le j <M} \frac{b_j}{j^2}  \sin  \frac{j\pi x}{b},  \label{eq:vide_v}
\end{eqnarray}
where $a_0,a_{-1}$ are two constant and can be determined by boundary conditions as shown in Eq (\ref{eq:a01}) below. 

The following notations and conventions will be adopted in the rest of this paper. A $k$-dim vector is considered as $(k,1)$ dimensional matrix unless specified otherwise. Define
\begin{eqnarray*} 
	x_k &=& k \lambda, \quad \lambda = \frac{b}M ,\quad 0\le k \le M, \\
	u_k &=& \tilde{u}_M(x_k), \quad v_k = \tilde{v}_M(x_k),  \quad z_k =\tilde{z}_M(x_k), \\
	p_k &=& p(x_k)h(x_k), \quad q_k=q(x_k)h(x_k), \quad r_k = r(x_k)h(x_k), \\
	f_k &=& f_h(x_k, v_k, u_k),  \quad g_k =g_h(x_k, v(s:e)), \\
	\mu_k &=& \mu(x_k)h(x_k), \quad \quad {\nu} = (\mu_k)_{1\le k< M},\\
	X &=& (x_k)_{1\le k < M}, \quad K = (1,2,\cdots, M-1)^T, \\
	 U &=& (u_k)_{1\le k< M},\quad V = (v_k)_{1\le k < M}, \quad Z = (z_k)_{1\le k < M}, \\
	F &=& (f_k)_{1\le k < M}, \quad G = (g_k)_{1\le k < M} , \quad B=(b_i)_{1\le i<M}\\
	I &=& (1,1,\cdots, 1)^T_{M-1}, \quad I_a = (-1,1, -1,\cdots, -1)^T_{M-1}, \\
	P &=& (p_k)_{1\le k< M}, \quad Q = (q_k)_{1\le k< M}, \quad R= (r_k)_{1\le k< M}.
\end{eqnarray*}
For any two matrices $A,B$ with same shape,  $A\circ B$ denotes the Hadamard product, which applies the element-wise product to two matrices. $A\cdot B$ denote the standard matrix multiplication when applicable. $A(i,:)$ and $A(:,j)$ is used to denote the $i$-th row and $j$-th column of $A$ respectively.  $W(k:l)$ denote the vector $(w_k, \cdots, w_l)^T$. In addition, $diag (W)$ is the diagonal matrix constructed by $W$.  Note we have
\[
s = x_{m}, \qquad e=x_{m+n}.
\]
At grid points of interpolation, ODE dynamic (\ref{eq:nonlinear_ode_order2_F}) is characterized by
\begin{equation} \label{eq:bdp_orde2_d_nonlinear} 
	Z = F + G.  
\end{equation}
We outline the rest of algorithm development in the following steps.
\begin{enumerate} 
	\item Use the analytic expressions (\ref{eq:vide_z}-\ref{eq:vide_v}) of $\tilde z,\tilde u,\tilde v$ to represent $Z$ and $U$ as linear functions of $V$ as shown in Eqs. (\ref{eq:vide_z_by_v}) and (\ref{eq:UAV})
	\item Express the integral component $G$ as a linear function of $V$ as shown in Eq. (\ref{G_v_integrable}). This is done in Subsection \ref{subsec:fide} based on trigonometrical interpolation of the associated kernel function. 
	\item Construct the matrix $\Phi$ and the vector $\Psi$ such that Eq (\ref{eq:bdp_orde2_d_nonlinear}) is equivalent to $\Phi V = \Psi$  as outlined in Algorithm \ref{Alg:VIDE}.
\end{enumerate}
By Eq. (\ref{eq:vide_z}-\ref{eq:vide_v}),  we obtain $\{z_k, u_k, v_k\}_{0\le k \le M}$ in $B$.
\begin{eqnarray}
	z_k &=& \sum_{0\le j<M}b_j \sin \frac{2\pi jk}{N},   \label{eq:z_M_d}\\
	u_k &=& a_0 - \frac{b}{\pi }\sum_{1\le j <M} \frac{b_j}{j} \cos  \frac{2\pi jk }{N},  \label{eq:tildeu_d}\\
	v_k &=& a_{-1} + a_0 x_k-(\frac{b}{\pi })^2 \sum_{1\le j <M}\frac{b_j}{j^2}\sin \frac{2\pi jk }{N}. \label{eq:tildev_d}
\end{eqnarray}
Eq (\ref{eq:tildeu_d})-(\ref{eq:tildev_d}) can be used to solve $a_0$ and $a_{-1}$:
\begin{equation}\label{eq:a01}
	a_0 = \frac{v_M - v_0}b, \quad a_{-1} = v_0.
\end{equation}
Define 
\[
S =(\sin\frac{2\pi jk}{N})_{1\le j,k <M}, \qquad C =  (\cos\frac{2\pi jk}{N})_{1\le j,k < M},
\]
and 
\[ 
\quad  O=\sqrt{\frac 2M}S,\qquad \Theta = O diag(1/K^2) O := (\theta_{ij})_{1\le i, j < M} .
\]
$B$ can be solved by $V$ as
\begin{equation}\label{eq:B}
	B = diag(K^2) S (\frac{2a_{-1} \pi^2}{Mb^2} I + \frac{2a_0 \pi^2}{bM^2}K - \frac{2\pi^2}{Mb^2} V).
\end{equation}
and therefore
\begin{equation}\label{eq:vide_z_by_v}
Z = SB=\Theta^{-1}(\frac{a_{-1} \pi^2}{b^2} I + \frac{a_0 \pi^2}{bM}K - \frac{\pi^2}{b^2} V)
\end{equation}
Let $U_e=(u_0,\dots, u_M)^T$, $V_e=(v_0,\dots, v_M)^T$, it is shown in  \cite{zou_tri_IV} that $U_e$ is covered by $V_e$ through linear transform   
\begin{equation}\label{eq:UAV}
	U_e = AV_e.
\end{equation}
where matrix $A$ is calculated by
\begin{eqnarray*}
	& & a_{0,0} = \frac{\pi}{b} S(I_a \circ \cot(\pi K/N))-\frac{\pi}{bM} S(I_a\circ K\circ \cot(\pi K/N))-\frac1b \label{eq:top}\\
	& & a_{0,1:M-1} = -\frac{\pi}{b} I'_a \circ \cot(\pi K'/N) \nonumber \\
	& & a_{0,M} = \frac{\pi}{bM} S(I_a\circ K\circ \cot(\pi K/N))+\frac1b \nonumber\\
	& & a_{i,0} = \frac{\pi}{2b} S((-1)^i \cot(i,:)I_a)- \frac{\pi}{2bM} S( (-1)^i I_a \circ \cot(i,:) \circ K) -\frac1b, \label{eq:middle}\\
	& & a_{i,1:M-1} = \frac{\pi}{2b}(-1)^{i+1} I_a' \circ \cot(i,:), \nonumber \\
	& & a_{i,M} = \frac{\pi}{2bM} S((-1)^i I_a \circ \cot(i,:) \circ K)+1/b. \nonumber \\
	& & a_{M,0} = -\frac{\pi}{b} S(I_a \circ \tan(\pi K/N))-\frac{\pi}{bM} S((K\circ \cot(\pi K/N))-\frac1b, \label{eq:bottom}\\
	& & a_{M,1:M-1} = \frac{\pi}{b} I'_a \circ \tan(\pi K'/N), \nonumber\\
	& & a_{M,M} = \frac{\pi}{bM} S(K\circ \cot(\pi K/N))+\frac1b,  \nonumber
\end{eqnarray*}
where $0<i<M$,  $S(V)$ denotes the summation over elements in a vector $V$,
\[
\cot(k, i) := Cot\frac{k+i}{N}\pi + Cot\frac{k-i}{N}\pi,
\]
and $Cot(x)=\cot(x)$ if $x/\pi$ is not integer and $Cot(x)=0 $ otherwise. By Eq. (\ref{eq:vide_z_by_v}), Eq (\ref{eq:bdp_orde2_d_nonlinear}) is equivalent to 
\begin{equation}\label{eq:ode_discrete_nonlinear}
	\frac{v_0\pi^2}{Mb^2} (MI-K) + \frac{v_M\pi^2}{M b^2} K -\frac{\pi^2}{b^2} V = \Theta \cdot (R + Q\circ V  + P\circ U +  G).
\end{equation}
With Eq (\ref{eq:UAV}), Eq. (\ref{eq:ode_discrete_nonlinear}) represents a linear algebraic system in $V_e$ if $G$  can be also covered by $V_e$ through linear transformation, which is to be developed in Subsection \ref{subsec:fide} and \ref{subsec:fide_en}. 
\subsection{Integral component with continuous kernel}\label{subsec:fide} 
We interpolate $k(x,y)$ by $2$-dim $sin$ polynomial as shown in \cite{zou_tri_FIDE}
\[
	K(x,t) = \sum_{1\le i,j<M} c_{ij} \sin \frac{\pi i x}b\sin \frac{\pi j t}b,
\]
and apply Eq (\ref{eq:vide_v}) to estimate $g_k$ for $1\le k <M$
\begin{eqnarray}
	\frac{g_k}{\mu_k}&=& \sum_{1 \le j<M} \frac {b \hat\eta(k,j)}{\pi }  (w_{-1}(k,j)a_{-1} + w_0(k,j) a_0) \nonumber \\
	&-&\sum_{1\le j,l<M}  \frac {b^3 \eta(k,j)}{2\pi^3} b_l w^k(j,l)    \label{v_g_continuous}
\end{eqnarray}
where
\begin{eqnarray*}
	\hat\eta_{kj} &=& \sum_{i}c_{ij}/j\sin\frac{2\pi ik}{N}, \quad \eta_{kj} = \sum_{i}c_{ij} \sin\frac{2\pi ik}{N},  \\
	w_{kj,-1} &=& \cos\frac{2\pi j m}{N} -  \cos\frac{2\pi j k}{N} \label{v_w_1},\\
	w_{kj,0} &=& s\cos\frac{2\pi j m}{N} -  x_k\cos\frac{2\pi j k}{N} + \frac{b}{\pi j} (\sin\frac{2\pi j k}{N} -\sin\frac{2\pi j m}{N}  )\label{v_w_0},\\
	w^k_{jl}&=& \frac{\delta_{j\neq l}}{(j-l)l^2} (\sin\frac{2\pi (j-l) k}{N} -\sin\frac{2\pi (j-l) m}{N} )  \\
	&+&  \frac{2\pi n}{l^2 \pi}\delta_{j=l} -\frac{1}{ (j+l)l^2} (\sin\frac{ (j+l) k}{N} -\sin\frac{2\pi (j+l) m}{N} ). 
\end{eqnarray*}
Define
\begin{eqnarray*}
	\eta &=& (\eta_{kj})_{(M-1, M-1)},  \quad \hat \eta = (\hat \eta_{kj})_{(M-1, M-1)}, \\
	W_{-1} &=& (w_{kj,-1})_{(0< k, j<M)}, \quad W_{0}=(w_{kj,0})_{(0<j, k<M)}.
\end{eqnarray*}
and for $1\le k<M$,
\[
W^k = (w^k_{jl})_{(0< j, l<M)}, \quad  H(k,:)=\eta(k,:)\cdot W^k\cdot diag(K^2)S.
\]
For any matrix $A$, let $A_s$ denote the column vector constructed by the summation of all columns of $A$. By Eq. (\ref{eq:a01}) and (\ref{eq:B}),  Eq (\ref{v_g_continuous}) can be linearly represented by $v_0$, $v_M$ and $V$ as   
\begin{equation}\label{G_v_continuous}
	G=  C_{0} v_0 + C_M v_M +  C_{rest} V,
\end{equation}
where
\begin{eqnarray*}
	C_M &=& \frac{1}{\pi} \nu\circ  (\hat \eta \circ W_0)_s -\frac{b}{\pi M^2} \nu \circ (H\cdot K),\\
	C_0 &=& \frac{b}{\pi}  \nu\circ (\hat \eta \circ W_{-1})_s - C_M  -\frac{b}{M\pi} \nu\circ (H\cdot I),\\
	C_{rest} &=& \frac{b}{\pi M} diag(\nu)\cdot H.
\end{eqnarray*}
\subsection{Integral component with integrable kernel}\label{subsec:fide_en} 
We assume that there are some singularity such that $k(x,y)=(|x-y|^{\gamma})\kappa(x,y)$ for $\gamma>-1$ and there are $k_1(x,y)$ and $k_2(x,y)$ such that $k_1(x,x)=0$ and $k_2(x,x)=0$ and 
\[
\frac{\partial {k_1(x,y)}}{\partial y}= k(x,y), \quad \frac{\partial {k_2(x,y)}}{\partial y}= k_1(x,y)
\]
As example for $k(x,y)=|x-y|^{\gamma}$, we have
\begin{eqnarray*}
	k_1(x,y)=\left\{\begin{array}{cc}
		-\frac{|x-y|^{1+\gamma}}{1+\gamma} & y\le x, \\
		\frac{|x-y|^{1+\gamma}}{1+\gamma} &  y> x\\
	\end{array}\right.  
\end{eqnarray*}
and 
\[
k_2(x,y) = \frac{1}{(1+\gamma)(2+\gamma)}|x-y|^{(1+\gamma)(2+\gamma)}.
\]
So we have
\[
	g(x,v(s:x)) = -v(s)k_1(x,s)+u(s)k_2(x,s) + \mu(x) \int^x_s k_2(x,t)z(t).
\]
We extend $k_2(x,y)$ to $[-b,b]^2$ by $2$-dim $sin$ polynomial    
\[
	k_2(x,t) = \sum_{1\le i,j<M} c_{ij} \sin \frac{\pi i x}b\sin \frac{\pi j t}b,
\]
and apply Eq (\ref{eq:vide_v}) to estimate $g_k$ for $1\le k <M$: 
\[
	\frac{g_k}{\nu_k}= -k_1(x_k, s) v_m  + k_2(x_k,s)u_m + \frac {b}{2\pi} \sum_{1\le j,l<M} \eta_{kj} w^k_{jl}b_l,    
\]
where
\begin{eqnarray*}
	\eta_{kj} &=& \sum_{i}c_{ij} \sin\frac{2\pi ik}{N},  \\
	w^k_{jl}&=& \frac{\delta_{j\neq l}}{(j-l)} (\sin\frac{2\pi (j-l) k}{N} -\sin\frac{2\pi (j-l) m}{N} ) \\
	&+&\frac{\pi}{b}(x_k-s)  -\frac{1}{ (j+l)} (\sin\frac{ 2\pi(j+l) k}{N} -\sin\frac{2\pi (j+l) m}{N} )\label{v_w}
\end{eqnarray*}
Define 
\begin{eqnarray*}
	\eta &=& (\eta_{kj})_{0< k, j<M},  	\quad W^k = (w^k_{jl})_{0< j, l<M}, \quad \\
	H(k,:)&=&\eta(k,:)\cdot W^k\cdot diag(K^2)S, \quad 1\le k<M.
\end{eqnarray*}
We have   
\begin{equation}\label{G_v_integrable}
	G=  C_{0} v_0 + C_mv_m + C_M v_M +  C_{rest} V 
\end{equation}
where
\begin{eqnarray*}
	& & C_0 = \frac{\pi}{b} \nu \circ (HI/M - HK/M^2) +  A(m,0)k_2(X,s)\\
	& & C_M = \frac{\pi}{bM^2} \nu \circ (HK) + A(m, M)k_2(X,s)\\
	& & C_m = -\nu\circ k_1(X,s), \\
	& & C_{rest} = -\frac{\pi}{b M} diag(\nu)\cdot H  + k_2(X,s)A(m, 1:M-1)\nonumber
\end{eqnarray*}
\subsection{Algorithm of solving VIDE}
\label{subsec:transform_linear_algebraic_system}
With linear representation of $G$ either by Eq (\ref{G_v_continuous}) or Eq (\ref{G_v_integrable}), we can solve VIDE (\ref{eq:vide_ode}-\ref{eq:vide_bdy2}) by Algorithm \ref{Alg:VIDE}.
\begin{Alg}\label{Alg:VIDE} Trigonometric Interpolation Based Approach for VIDE. 
	\begin{enumerate}
		\item Construct the following $M+1$ dimensional linear algebraic system  that consists of Eq (\ref{eq:nonlinear_ode_order2:diri_F}) (first in the system), Eq (\ref{eq:ode_discrete_nonlinear}) (last in the system), and Eq (\ref{G_v_continuous}) or Eq (\ref{G_v_integrable}). 
		\begin{equation}\label{eq:PhiVG}
			\Phi V_e = \Psi,  \quad \Phi = (\phi_{jk})_{0\le j,k\le M}, \Psi = (\psi_j)_{0\le j\le M}.
		\end{equation}
		\item  Calculate $\Psi=(\psi_i)_{0\le i\le M}$ by 
		\begin{equation} \label{eq:psi}
			\psi_0=\alpha, \qquad \psi_M= \beta, \quad  \Psi(1:M-1)= -\Theta R.
		\end{equation}
		\item Calculate $\Phi=(\phi_{ij})_{0\le i,k \le M}$ in three steps
		\begin{enumerate}
			\item Identify the first and last row of $\Phi$, which are associated to Eq  (\ref{eq:nonlinear_ode_order2:diri_F}) and Eq (\ref{eq:ode_discrete_nonlinear}) respectively  
			\begin{eqnarray*}
				& & \phi_{0,k} = d_{12}\cdot A(m,k) + d_{14}\cdot A(m+n,k) + \delta_{m,k}d_{11} + \delta_{m+n,k}d_{13}, \label{eq:phi_0} \\
				& &\phi_{M,k} = d_{22}\cdot A(m,k) + d_{24}\cdot A(m+n,k) + \delta_{m,k}d_{21} + \delta_{m+n,k}d_{23}. \label{eq:phi_M} 	 
			\end{eqnarray*}
			\item Identify other $M-1$ rows of $\Phi$. For $0<i<M$, 
			\begin{eqnarray*}
			& &	\phi_{i,0} = -\frac{(M-i)\pi^2}{Mb^2} + (\theta(i,:)\circ P^T) \cdot A(1:M-1,0) + \theta(i,:) C_0,\label{eq:phi_i_0} \\
			& &	\phi_{i,M} = -\frac{i \pi^2}{b^2 M} + (\theta(i,:) \circ P^T) \cdot A(1:M-1, M) + \theta(i,:)C_M,\label{eq:phi_i_M} \\
			& &	\phi(i,1:M-1) = \frac{\pi^2}{b^2}(\delta_{ij})_{1\le j<M} +\theta(i,:)\circ Q^T  \\
				&+&(\theta(i,:) \circ P^T) \cdot A(1:M-1,1:M-1) \label{Phi}  + \theta(i,:)C_{rest}.
			\end{eqnarray*}
			\item  Add the impact of $C_m$ and $C_{m+n}$ if Eq (\ref{G_v_integrable}) is used to estimate $G$ component.
			\begin{eqnarray*}
			& &	\Phi(1:M-1,m) \leftarrow \Phi(1:M-1,m) + \theta\cdot c_m\\
			& &	\Phi(1:M-1,m+n) \leftarrow \Phi(1:M-1,m+n) + \theta\cdot c_{m+n}\\
			\end{eqnarray*}
		\end{enumerate}
		\item Solve Eq (\ref{eq:PhiVG}) with calculated $\Phi$ and $\Psi$ in previous steps.
		\item Apply Eq (\ref{eq:a01}) and Eq (\ref{eq:B}) to compute coefficients in Eq (\ref{eq:vide_v}) for the approximation $v_M$ on solution of VIDE (\ref{eq:vide_ode}-\ref{eq:vide_bdy2}).
	\end{enumerate}
\end{Alg}
\newpage
\section{Numerical examples}\label{sec:numerical_examples}
In this section,  we  test the method proposed in Section \ref{sec:fide} with kernel functions $k(x,y)$ defined in Table \ref{tab:labels_test_dim2} with focus on $K_1$. 
\begin{table}[htbp]
\caption{The selected kernels for the performance test on Algorithm \ref{Alg:VIDE}                   }
\begin{tabular}{cccc}
	Notation & $k(x,y)$ & Notation & $k(x,y)$ \\ \hline \hline
	$K_1$  & $|x-y|^\gamma$ & $Exp$ & $exp(x+y)$   \\ \hline
	$K_2$  & $|x^2-y^2|^\gamma$ & $Sin$   & $sin(x+y)$  \\ \hline
\end{tabular}%
\label{tab:labels_test_dim2}%
\end{table}%
For any given smooth functions $p(x),q(x),\mu(x)$, matrix $(d_{ij})_{2\times 4}$ and $f(x)\in C^2([s, e])$, one can see that $f(x)$ solves the following VIDE
\begin{eqnarray*}
	y'' &=& p(x) y' +  q(x)y  + r(x) + \mu(x) \int^e_sk(x,t)y(t)dt\\
	\alpha  &=& d_{11}y(s) + d_{12}y'(s) + d_{13}y(e) + d_{14}y'(e) \\
	\beta  &=& d_{21}y(s) + d_{22}y'(s) + d_{23}y(e) + d_{24}y'(e) 
\end{eqnarray*}
where $r(x)$ and constant $\alpha,\beta$ are determined by
\begin{eqnarray*}
	r(x) &=& f''(x)  - p(x)f'(x) - q(x)f(x) -  \mu(x)\int^e_s k(x,t)f(t)dt, \\
	\alpha  &=& d_{11}f(s) + d_{12}f'(s) + d_{13}f(e) + d_{14}f'(e), \\
	\beta  &=& d_{21}f(s) + d_{22}f'(s) + d_{23}f(e) + d_{24}f'(e). 
\end{eqnarray*}
All tests are based on the following parameters 
\[
s=1,\quad e=3, \quad \delta=1, \quad p(x) \equiv  0.1, \quad q(x) \equiv 1, \quad \mu(x) \equiv 1,
\]
target function $f$ can be one of following list 
\begin{equation}\label{test_functions}
	\cos (\frac{\pi x}2 ) , \quad  \cos (\frac{3\pi x}2) , \quad e^x, \quad x^2. 
\end{equation}
We consider four sets of boundary conditions in Table \ref{tab:test}.
\begin{table}[htbp]
	\centering
	\small
	\caption{The types of boundary conditions. $\{d_{ij}\}_{1\le i,j\le 4}$ are parameters in Eq (\ref{eq:vide_bdy1}-\ref{eq:vide_bdy2}).}
	\begin{tabular}{lrrrrrrrrl}
		type	& \multicolumn{1}{l}{$d_{11}$} & \multicolumn{1}{l}{$d_{12}$} & \multicolumn{1}{l}{$d_{13}$} & \multicolumn{1}{l}{$d_{14}$} & \multicolumn{1}{l}{$d_{21}$} & \multicolumn{1}{l}{$d_{22}$} & \multicolumn{1}{l}{$d_{23}$} & \multicolumn{1}{l}{$d_{24}$} & condition on \\ \hline\hline
		$Neumann$ & 1     & 0     & 0     & 0     & 0     & 1     & 0     & 0     & $v_s,u_s$ \\
		$Dirichlet$ & 1     & 0     & 0     & 0     & 0     & 0     & 1     & 0     & $v_s, v_e$ \\
		$Mix_1$   & 1     & 0     & 0     & 0     & 0     & 0     & 0     & 1     & $v_s, u_e$ \\ 
		$Mix_2$ & 1     & 1     & 0     & 0     & 0     & 0     & 1     & 1     & $v_s+u_s, v_e+u_e$ \\ \hline
	\end{tabular}%
	\label{tab:test}%
\end{table}%
The performance is measured by the normalized max error $max_{e}$ defined as
\[
err_{e} = \frac{\max_{x\in S_e}|f(x)-\tilde{v}_M(x)|}{\max_{x\in S_e}|f(x)|}, \quad  S_e = \{ x^p_k, 0\le k \le 2^{11}  \} \cap [s,e]
\]  
where $S_e=\{x^p_k\}$ consists equally spaced points over $[-b,b]$ based on $q=10$ and contains more than interpolation grid points used to determine $\tilde v_M$ since $M$ is always less than $2^{10}$ in all test cases. 
\subsection{Numerical results on performance of continuous kernels}\label{subsec:numerical_results_continous}
The performance of four continuous kernels in Table \ref{tab:labels_test_dim2} are tested with target function $f(x)=\cos(3\pi x/2)$ and Dirichlet boundary condition,  and the results are shown in Table \ref{tab:volterra_en_convergenc}.  Almost same accuracy is observed across four kernel functions. Note that trigonometric interpolation error for $K_1$ and $K_2$ is much higher than the other two kernels. The testing result might imply that the accuracy of VIDE's solution is mainly determined by the performance of trigonometric interpolation of kernel at grid points, where four kernels exhibits similar performance. Performance of trigonometric interpolation for the kernels in Table \ref{tab:labels_test_dim2} can be found in \cite{zou_tri_FIDE}.
\begin{table}[htbp]
\caption{The performance of covered kernels with   $f(x)=\cos (\frac{3\pi }2 x)$}
\begin{tabular}{lrrrrr}
	kernel & type  & \multicolumn{1}{l}{$\gamma$} & \multicolumn{1}{l}{$max_{e}$} &  \multicolumn{1}{l}{q} \\ \hline\hline
	$K_1$ & Dirichlet & 0.5   & 5.0E-11  & 7 \\ \hline
	$K_2$ & Dirichlet & 0.5  & 5.0E-11 &  7 \\ \hline
	$Exp$ & Dirichlet & NA   & 5.0E-11  & 7 \\ \hline
	$Sin$ & Dirichlet & NA   & 5.0E-11  & 7 \\ \hline
\end{tabular}%
\label{tab:volterra_en_convergenc}%
\end{table}%

Table \ref{tab:test_on_boundary} shows the performance of different boundary conditions with target function $f(x)=\cos(3\pi x/2)$ for VIDE with kernel $K_1$. The performance is decent for all cases and reaches the best accuracy with Dirichlet condition. 
\begin{table}[htbp]
\caption{The performance of covered boundary conditions with   $f(x)=\cos (\frac{3\pi }2 x)$}
\begin{tabular}{lrrrrr}
	kernel & type  & \multicolumn{1}{l}{$\gamma$} & \multicolumn{1}{l}{$max_{e}$} & \multicolumn{1}{l}{q} \\ \hline\hline
	$K_1$ & Neumann & 0.5    & 2.2E-08 &  7 \\ \hline
	$K_1$ & Dirichlet & 0.5  & 5.0E-11 &  7 \\ \hline
	$K_1$ & $Mix_1$ & 0.5    & 3.9E-09 &  7 \\ \hline
	$K_1$ & $Mix_2$ & 0.5    & 3.4E-08 &  7 \\ \hline 
\end{tabular}%
\label{tab:test_on_boundary}%
\end{table}%

Table \ref{tab:test_on_convergence} shows convergence property and high accuracy is achieved for $q\ge 7$ with Dirichlet condition. As expected,  the performance is improved as $q$ increases. 
\begin{table}[htbp]
\caption{Convergence Test with $f(x)=\cos (\frac{3\pi }2 x)$}
\begin{tabular}{lrrrr}
	kernel & type  & \multicolumn{1}{l}{$\gamma$} & \multicolumn{1}{l}{$max_{e}$} & \multicolumn{1}{l}{q} \\ \hline\hline
	$K_1$ & Dirichlet & 0.5   & 4.8E-03 & 4 \\ \hline
	$K_1$ & Dirichlet & 0.5   & 9.2E-05 & 5 \\ \hline
	$K_1$ & Dirichlet & 0.5   & 1.6E-07 & 6 \\ \hline
	$K_1$ & Dirichlet & 0.5   & 5.0E-11 & 7 \\ \hline
	$K_1$ & Dirichlet & 0.5   & 5.5E-14 & 8  \\ \hline
\end{tabular}%
\label{tab:test_on_convergence}%
\end{table}%
\subsection{Numerical results on performance of kernel with singularities}\label{subsec:numerical_results_sigularity}
This subsection includes the similar testing results as in Subsection \ref{subsec:numerical_results_continous} for VIDEs with kernel $K_1$ that bears singularities when $\gamma<0$. 

Table \ref{tab:volterra_en_gamma} shows the impact of degree of singularity of $k(x,y)$. Similar performance is observed for all covered $\gamma$, consistent to what has been observed in Subsection \ref{subsec:numerical_results_continous}.
\begin{table}[htbp]
\caption{Impact of kernel's smoothness with $\mu=1,\theta=3\pi/2, M=2^q$}
\begin{tabular}{lrrrr}
	kernel & type & \multicolumn{1}{l}{$\gamma$} & \multicolumn{1}{l}{$max_{e}$} & \multicolumn{1}{l}{$q$} \\ \hline\hline
	$K_1$ & Dirichlet & -0.9   & 4.5E-08 & 7 \\\hline
	$K_1$ & Dirichlet & -0.5   & 7.3E-09 & 7 \\\hline
	$K_1$ & Dirichlet & 0     & 7.5E-09 & 7 \\\hline
	$K_1$ & Dirichlet & 0.5  & 8.0E-09 & 7 \\\hline
	$K_1$ & Dirichlet & 1.5  & 9.8E-09 & 7 \\\hline
	$K_1$ & Dirichlet & 2    & 1.1E-08 & 7 \\ \hline
\end{tabular}%
\label{tab:volterra_en_gamma}%
\end{table}%
Table \ref{tab:volterra_en_convergence} shows convergence property and high accuracy is achieved for $q\ge 7$ with Dirichlet condition. As expected, the performance is improved as $q$ increases. 
\begin{table}[htbp]
\caption{Convergence Test with $f(x)=\cos (\frac{3\pi }2 x)$}
\begin{tabular}{lrrrr}
	kernel & type & \multicolumn{1}{l}{$\gamma$} & \multicolumn{1}{l}{$max_{e}$} & \multicolumn{1}{l}{$q$} \\ \hline\hline
	$K_1$ & Dirichlet & -0.5   & 4.5E-03 & 4 \\\hline
	$K_1$ & Dirichlet & -0.5   & 8.8E-05 & 5 \\\hline
	$K_1$ & Dirichlet & -0.5   & 1.7E-07 & 6 \\\hline
	$K_1$ & Dirichlet & -0.5   & 7.3E-09 & 7 \\\hline
	$K_1$ & Dirichlet & -0.5   & 4.4E-10 & 8 \\\hline
	$K_1$ & Dirichlet & -0.5   & 2.5E-11 & 9 \\ \hline
\end{tabular}%
\label{tab:volterra_en_convergence}%
\end{table}%
Table \ref{tab:volterra_en_boundary} shows the performance of different boundary conditions with target function $f(x)=\cos(3\pi x/2)$ for VIDE with kernel $K_1, \gamma=-0.5$. The performance is decent for all cases and reaches the best accuracy with Dirichlet condition, consistent to what is observed in Subsection \ref{subsec:numerical_results_continous}.
\begin{table}[htbp]
\caption{The performance of covered boundary conditions with   $f(x)=\cos (\frac{3\pi }2 x)$}
\begin{tabular}{llrrr}
	kernel & type & \multicolumn{1}{l}{$\gamma$} & \multicolumn{1}{l}{$max_{e}$} & \multicolumn{1}{l}{$q$} \\ \hline\hline
	$K_1$ & Neumann & -0.5   & 4.1E-07 & 7 \\\hline
	$K_1$ & Dirichlet & -0.5   & 7.3E-09 & 7 \\\hline
	$K_1$ & $Mix_1$ & -0.5   & 7.5E-09 & 7 \\\hline
	$K_1$ & $Mix_2$ & -0.5   & 1.2E-07 & 7 \\ \hline
\end{tabular}%
\label{tab:volterra_en_boundary}%
\end{table}%
Table \ref{tab:volterra_en_boundary} provides the performance of VIDE with kernel $K_1, \gamma=-0.5$ with covered target functions.
High accuracy is reached for all cases. 
\begin{table}[htbp]
\caption{Performance with covered target functions}
\begin{tabular}{llrrrr}
	kernel & type & \multicolumn{1}{l}{$\gamma$} &  target function & \multicolumn{1}{l}{$max_{e}$}  & \multicolumn{1}{l}{$q$}  \\ \hline\hline
	$K_1$ & Dirichlet & -0.5  & $cos(\pi x/2)$    & 2.6e-09  & 7 \\\hline
	$K_1$ & Dirichlet & -0.5  & $cos(3\pi x/2)$    &  7.3e-09  & 7 \\\hline
	$K_1$ & Dirichlet & -0.5  & $x^2$    & 1.9e-09  & 7 \\\hline
	$K_1$ & Dirichlet & -0.5  & $exp(x)$    & 4.0e-9 & 7 \\  \hline
\end{tabular}%
\label{tab:volterra_en_multiple}%
\end{table}%

\section{Summary}\label{sec:summary}
In this paper, we proposed an trigonometric interpolation method, Algorithm \ref{Alg:VIDE}, for the numerical solution of a second order VIDE.  The new method converts VIDE to a linear algebraic system and it bears several advantages as described in Section \ref{sec:vide_intro}.  The tests on Algorithm \ref{Alg:VIDE} in Section \ref{sec:numerical_examples} cover various boundary conditions with both continuous and integrable kernels. Decent performance is observed across all covered scenarios with a moderate size of grid points. 

\end{document}